\documentclass[11pt,reqno]{amsart}
\usepackage{enumerate}
\usepackage{amsfonts, amsmath, amssymb, amscd, amsthm, bm}
\usepackage{enumerate}
\usepackage{url}
\usepackage[linktocpage=true,colorlinks,citecolor=magenta,linkcolor=blue,urlcolor=magenta]{hyperref}
\usepackage{multicol}
\usepackage{comment}

 \newtheorem{thm}{Theorem}[section]

 \newtheorem{lem}[thm]{Lemma}
 \newtheorem{assump}[thm]{Assumption}
 \newtheorem{prop}[thm]{Proposition}
 \theoremstyle{definition}
 
 \theoremstyle{remark}
 \newtheorem{rem}[thm]{Remark}

 \theoremstyle{claim}
 \numberwithin{equation}{section}
 \newcommand{\RR}{{\mathbb R}}

\begin{document}
\title{Quermassintegral preserving curvature flow in Hyperbolic space}
\author[B. Andrews]{Ben Andrews}
\author[Y. Wei]{Yong Wei}
\address{Mathematical Sciences Institute,
Australian National University,
ACT 2601 Australia}
\email{\href{mailto:Ben.andrews@anu.edu.au}{Ben.andrews@anu.edu.au}, \href{mailto:yong.wei@anu.edu.au}{yong.wei@anu.edu.au}}

\date{\today}
\subjclass[2010]{53C44; 53C21}
\keywords {Quermassintegral preserving flow, hyperbolic space, Alexandrov reflection.}
\thanks{The research was supported by Australian Laureate Fellowship FL150100126 of the Australian Research Council.}


\begin{abstract}
We consider the quermassintegral preserving flow of closed \emph{h-convex} hypersurfaces in hyperbolic space with the speed given by any positive power of a smooth symmetric, strictly increasing, and homogeneous of degree one function $f$ of the principal curvatures which is inverse concave and has dual $f_*$ approaching zero on the boundary of the positive cone. We prove that if the initial hypersurface is \emph{h-convex}, then the solution of the flow becomes strictly \emph{h-convex} for $t>0$, the flow exists for all time and converges to a geodesic sphere exponentially in the smooth topology.
\end{abstract}

\maketitle

\section{Introduction}
Let $X_0: M^n\to \mathbb{H}^{n+1}$ be a smooth embedding such that $M_0=X_0(M)$ is a closed smooth hypersurface in the hyperbolic space $\mathbb{H}^{n+1}$. We consider the smooth family of immersions $X:M^n\times [0,T)\rightarrow \mathbb{H}^{n+1}$ satisfying
\begin{equation}\label{flow-VMCF}
 \left\{\begin{aligned}
 \frac{\partial}{\partial t}X(x,t)=&~(\phi(t)-\Psi(\mathcal{W}(x,t)))\nu(x,t),\\
 X(\cdot,0)=&~X_0(\cdot),
  \end{aligned}\right.
 \end{equation}
 where $\nu(x,t)$ is the unit outward normal of $M_t=X(M,t)$, and $\Psi(\mathcal{W})=F^{\alpha}(\mathcal{W})$ where $F$ is a smooth invariant function of the Weingarten matrix $\mathcal{W}=(h_i^j)$ of $M_t$.  The global term $\phi(t)$ is chosen to keep one of the Quermassintegrals of the hypersurface constant (we will explain this below).  We assume that $\alpha>0$ and $F$ satisfies the following conditions:
\begin{assump}\label{s1:assum}
\begin{itemize}
\item[(i)] $F(\mathcal{W})=f(\kappa(\mathcal{W}))$, where $\kappa(\mathcal{W})$ gives the eigenvalues of $\mathcal{W}$ and $f$ is a smooth symmetric function on
    \begin{equation*}
      \Gamma_+=\{x=(x_1,\cdots,x_n)\in\mathbb{R}^n:\ ~x_i>0,~i=1,\cdots,n\}.
    \end{equation*}
\item[(ii)] $f$ is strictly increasing, i.e., $\dot{f}^i={\partial f}/{\partial \kappa_i}>0$ on $\Gamma_+$, $\forall~i=1,\cdots,n$.
  \item[(iii)] $f$ is homogeneous of degree $1$, i.e., $f(k\kappa)=kf(\kappa)$ for any $k>0$.
  \item[(iv)] $f$ is strictly positive on $\Gamma_+$ and is normalized such that $f(1,\cdots,1)=1$.
  \item[(v)] $f$ is inverse concave, i.e., the function
    \begin{equation}\label{s3:f-dual}
      f_*(x_1,\cdots,x_n)=f(x_1^{-1},\cdots,x_n^{-1})^{-1}
    \end{equation}
    is concave.
  \item[(vi)] $f_*$ approaches zero on the boundary of $\Gamma_+$.
 \end{itemize}
\end{assump}
To describe the global term $\phi(t)$ in \eqref{flow-VMCF}, we first recall the (normalized) $k$-th mean curvature $E_k$ of a smooth closed hypersurface $M$ and the quermassintegrals $W_k(\Omega)$ of the bounded domain $\Omega$ enclosed by $M$:  If $\Omega$ is a (geodesically) convex domain in $\mathbb{H}^{n+1}$, then the quermassintegrals of $\Omega$ are defined as follows (see \cite{Sant2004,Sol2006,WX}):
\begin{equation}\label{s2:quer-def}
  W_k(\Omega)=~\frac{(n+1-k)\omega_{k-1}\cdots\omega_0}{(n+1)\omega_{n-1}\cdots\omega_{n-k}}\int_{\mathcal{L}_k}\chi(L_k\cap\Omega)dL_k,\quad k=1,\cdots,n,
\end{equation}
where $\mathcal{L}_k$  is the space of $k$-dimensional affine subspaces $L_k$ in $\mathbb{H}^{n+1}$. The function $\chi$ is defined to be $1$ if  $L_k\cap\Omega\neq \emptyset$ and to be $0$ otherwise. In particular, we have
\begin{equation*}
  W_0(\Omega)=|\Omega|,\quad W_{n+1}(\Omega)=\frac{\omega_n}{n+1},\quad W_1(\Omega)=\frac 1{n+1}|\partial\Omega|.
\end{equation*}
If the boundary $M=\partial \Omega$ is smooth (at least of class $C^2$), we can define the principal curvatures $\kappa=(\kappa_1,\cdots,\kappa_n)$ as the eigenvalues of the Weingarten matrix $\mathcal{W}$ of $M$. For each $k\in\{1,\cdots,n\}$ the $k$-th mean curvature $E_k$ of $M$ is then defined as the normalized $k$-th elementary symmetric functions of the principal curvatures of $M$:
 \begin{equation*}
   E_{k}={\binom{n}{k}}^{-1}\sum_{1\leq i_1<\cdots<i_k\leq n}\kappa_{i_1}\cdots \kappa_{i_k}.
 \end{equation*}
These include $E_1=H/n=(\kappa_1+\cdots+\kappa_n)/n$ (the normalized mean curvature) and $E_n=\kappa_1\cdots\kappa_n$ (the Gauss curvature). The \emph{curvature integrals} of $\Omega$ are then defined by
\begin{equation}\label{s2:CurInt}
  V_{n-k}(\Omega)~=~\int_{\partial\Omega}E_k,\quad k=0,1,\cdots,n.
\end{equation}
The quermassintegrals and curvature integrals of smooth convex domain $\Omega$ in $\mathbb{H}^{n+1}$ are related as follows:
\begin{align}
  V_{n-k}(\Omega) =& (n+1)\left(W_{k+1}(\Omega)+\frac k{n+2-k}W_{k-1}(\Omega)\right),\quad k=1,\cdots,n \label{s2:Quer-int1}\\
  V_n(\Omega) =&(n+1)W_1(\Omega)=|\partial\Omega|.
\end{align}

Besides the (geodesic) convexity, there is a stronger notion of convexity for regions in $\mathbb{H}^{n+1}$: \emph{horospherical convexity}. A domain $\Omega$ in $\mathbb{H}^{n+1}$ is called \emph{h-convex} (or, \emph{horospherically convex}) if at every boundary point $p\in M=\partial\Omega$ there is a horoball $\mathcal{H}$ of $\mathbb{H}^{n+1}$ which contains $\Omega$ and touches at $p$ (i.e. $p$ is a boundary point of ${\mathcal H}$). Recall that a horoball in $\mathbb{H}^{n+1}$ is the union of those geodesics balls which have centre on a given geodesic ray from $p$ and which have $p$ as a boundary point.
A smooth domain $\Omega$ is  \emph{h-convex} in $\mathbb{H}^{n+1}$  if and only if all the principal curvatures of $M=\partial\Omega$ are bounded from below by $1$, and $\Omega$ is strictly \emph{h-convex} if all the principal curvatures of $M=\partial\Omega$ are strictly bigger than $1$. In this paper, when we say a hypersurface $M$ is (strictly) \emph{h-convex}, we mean that the domain $\Omega$ enclosed by $M$ is (strictly) \emph{h-convex}.

Fix an integer $k=1,\cdots,n$.  If we define the function $\phi(t)$ in \eqref{flow-VMCF} by
\begin{equation}\label{s1:phit}
  \phi(t)=\frac{\int_{M_t}E_k\Psi d\mu_t}{\int_{M_t}E_kd\mu_t},
\end{equation}
then the quermassintegral $W_k(\Omega_t)$ of $\Omega_t$ remains constant along the flow \eqref{flow-VMCF} (see \S \ref{sec:pre}). We call flows of the form \eqref{flow-VMCF} with $\phi(t)$ given by \eqref{s1:phit} \emph{quermassintegral preserving curvature flows}. In particular, in the case $k=0$, these are \emph{volume preserving curvature flows}. The main result of this paper is the following:

\begin{thm}\label{thm1-1}
Let $k\in \{1,\cdots,n\}$ and $X_0: M^n\to \mathbb{H}^{n+1}$ be a smooth embedding such that $M_0=X_0(M)$ is a closed h-convex hypersurface in $\mathbb{H}^{n+1}$. Then for any power $\alpha>0$, the flow \eqref{flow-VMCF} with $F$ satisfying Assumption \eqref{s1:assum} and the global term $\phi(t)$ given by \eqref{s1:phit} has a smooth h-convex solution $M_t$ for all time $t\in [0,\infty)$, and $M_t$ converges smoothly and exponentially to a geodesic sphere of radius $r_{\infty}$ determined by $W_k(B_{r_{\infty}})=W_k(\Omega_0)$  as $t\to\infty$.
\end{thm}

\begin{rem}
Some important examples of functions satisfying the Assumption \ref{s1:assum} include: (i). $f=E_k^{1/k}$ for all $k=1,\cdots,n$; (ii) the power-means $H_r=(\frac 1n\sum_ix_i^r)^{1/r}$ for all $r>0$; and (iii) any function of the form $E_1^\sigma G^{1-\sigma}$ where $0<\sigma<1$ and $G$ is homogeneous degree one, normalised, increasing in each argument, and inverse-concave. See \S \ref{sec:pre} for further discussion.
\end{rem}

Constrained curvature flows have been studied extensively in recent years. In 1987, Huisken \cite{huis-87} studied the volume preserving mean curvature flow in Euclidean space $\mathbb{R}^{n+1}$, and proved that starting from any strictly convex hypersurface a solution exists for all time $t\in[0,\infty)$ and converges smoothly to a round sphere. In 2001, the first named author \cite{And2001-Aniso} studied volume preserving anisotropic mean curvature flows in $\mathbb{R}^{n+1}$ and obtained a similar result. Later, McCoy \cite{McC2003, McC2004, McC2005,Mcc2017} studied some mixed volume preserving curvature flow driven by homogeneous of degree one curvature functions. For higher homogeneity, by imposing a strong pinching assumption on the initial hypersurface, Cabezas-Rivas and Sinestrari \cite{Cab-Sin2010} proved convergence results for the flow \eqref{flow-VMCF} in $\mathbb{R}^{n+1}$ with $\Psi=E_k^{\alpha/k}$ where $k=1,\cdots,n$ and $\alpha>1$.  Using the monotonicity of the isoperimetric ratio, Sinestrari \cite{Sin-2015} proved a convergence result for the flow \eqref{flow-VMCF} with  $\Psi=H^{\alpha}$ in $\mathbb{R}^{n+1}$ and for any positive power $\alpha>0$. This was generalized in \cite{Be-Sin2016} for volume (and area) preserving non-homogeneous mean curvature flow in $\mathbb{R}^{n+1}$.  Very recently, the authors \cite{And-Wei2017} removed the pinching assumption in \cite{Cab-Sin2010} and proved that the flow \eqref{flow-VMCF} in $\mathbb{R}^{n+1}$ with $\Psi=E_k^{\alpha/k}$ for $k=1,\cdots,n$ and any $\alpha>0$, will deform any strictly convex hypersurface to a round sphere smoothly.

The volume preserving flow in hyperbolic space $\mathbb{H}^{n+1}$ was firstly studied by Cabezas-Rivas and Miquel \cite{Cab-Miq2007} in 2007. By imposing \emph{h-convexity} on the initial hypersurface, they proved that the flow \eqref{flow-VMCF} with $\Psi=H$ exists for all time and converges smoothly to a geodesic sphere. This result was generalized recently by Bertini and Pipoli \cite{Be-Pip2016} to volume preserving non-homogeneous mean curvature flow. In particular, their result includes the flow with velocity given by $\Psi=H^\alpha$ with $\alpha>0$. Some mixed volume preserving flows were considered in \cite{Mak2012,WX} with $\Psi$ given by some homogeneous of degree one curvature function $F$. By assuming \emph{h-convexity} and strong pinching on the initial hypersurface, Guo-Li-Wu \cite{GuoLW-2013} proved the convergence of the flow \eqref{flow-VMCF} with $\Psi=E_k^{\alpha/k}, k=1,\cdots,n$ and power $\alpha>1$ by following the same procedure as the Euclidean case in \cite{Cab-Sin2010}.

The paper is organized as follows: In section \ref{sec:pre}, we collect some preliminaries on hyperbolic geometry, the evolution equations along the flow \eqref{flow-VMCF}, and examples of the function satisfying the Assumption \ref{s1:assum}.  In section \ref{sec:h-conv}, we prove that the \emph{h-convexity} is preserved along the flow \eqref{flow-VMCF} for inverse concave $f$  and for all power $\alpha>0$. To show this, we will apply the tensor maximum principle proved by the first named author in \cite{And2007} (which generalized Hamilton's \cite{Ha1982} theorem). In section \ref{sec:ub-F}, the preservation of $W_k(\Omega_t)$ and the \emph{h-convexity} will be used to estimate the inner radius and outer radius of $\Omega_t$. Then we adapt Tso's \cite{Tso85} technique to prove a uniform upper bound on $F$. In section \ref{sec:LTE}, by the assumption that $f_*$ approaches zero on the boundary of $\Gamma_+$ and the upper bound on $F$, we derive a uniform upper bound on the principal curvatures. The $h$-convexity together with the boundedness of principal curvatures makes the evolution equation uniformly parabolic.
By projecting the flow solution into the unit ball in $\mathbb{R}^{n+1}$ and using the Gauss map parametrization, we write the flow \eqref{flow-VMCF} as a scalar parabolic partial differential equation for the support function which is concave with respect to the second spatial derivatives. Then the H\"{o}lder estimate of Krylov-Evans \cite{Krylov} and the parabolic Schauder theory \cite{lieberman1996} can be applied to derive the higher order derivative estimates of the solution $M_t$.  From this we conclude that the solution of \eqref{flow-VMCF} exists for all time $t\in[0,\infty)$.

In previous work, the convergence of solutions as $t\to\infty$ was deduced using either the monotonicity of curvature pinching ratios \cite{huis-87,McC2003,McC2004,McC2005,Mcc2017,Cab-Miq2007,Cab-Sin2010,GuoLW-2013,Mak2012} or of isoperimetric ratios \cite{And2001-Aniso, Sin-2015,Be-Sin2016, Be-Pip2016, And-Wei2017}.  In our situation for general $F$ and $\alpha$, neither of these arguments is available.  Instead,
we apply the Alexandrov reflection method to prove that the hypersurfaces approach a sphere, and a linearisation argument to prove exponential convergence.

\begin{rem}
In the case where $F=E_k^{1/k}$, we have as in \cite{And-Wei2017} that the quermassintegral $W_k(\Omega_t)$ is monotone decreasing along the volume preserving flow for any $\alpha>0$. This can also be used to deduce the smooth convergence to the geodesic sphere. We describe this alternative argument in the appendix.
\end{rem}

\begin{rem}
As in \cite{Be-Sin2016,Be-Pip2016}, a result similar to Theorem \ref{thm1-1} is also true for non-homogeneous constrained flows \eqref{flow-VMCF} with $\Psi=\psi(f)$, where $f$ is a function satisfying Assumption \ref{s1:assum}, and $\psi: [0,+\infty)\to \mathbb{R}$ is in $C^0([0,\infty))\cap C^2((0,\infty))$ and satisfies
\begin{itemize}
  \item[(i)] $\psi(r)>0$, $\psi'(r)>0$ for all $r>0$;
  \item[(ii)] $\lim_{r\to\infty}\psi(r)=\infty$;
  \item[(iii)] $\lim_{r\to\infty}\frac{\psi'(r)r^2}{\psi(r)}=\infty$;
  \item[(iv)] $\psi''(r)r+2\psi'(r)\geq 0$ for all $r>0$
\end{itemize}
In fact, the flow is parabolic due to item (i); item (iv) is used to show that the \emph{h-convexity} is preserved along the flow (see Remark \ref{s3:rem}); items (ii)--(iii) are used to estimate the upper bound of $f$. The remaining proof is similar.
\end{rem}

\section{Preliminaries}\label{sec:pre}
In this section, we collect some preliminary results concerning hyperbolic geometry, the properties of inverse concave functions, and examples of functions satisfying Assumption \ref{s1:assum}.  We also collect the evolution equations for several geometric quantities of the solution $M_t$ of the flow \eqref{flow-VMCF}.

\subsection{Hyperbolic geometry}

Let $M$ be a smooth closed hypersurface in $\mathbb{H}^{n+1}$. We denote by $g_{ij}, h_{ij}$ and $\nu$ the induced metric, the second fundamental form and unit outward normal vector of $M$. The Weingarten map is denoted by $\mathcal{W}=(h_i^j)$, where $h_i^j=h_{ik}g^{kj}$. The principal curvatures $\kappa=(\kappa_1,\cdots,\kappa_n)$ of $M$ are defined as the eigenvalues of $\mathcal{W}$. As mentioned before, $M$ is \emph{h-convex} if and only if $\kappa_i\geq 1$ for all $i=1,\cdots,n$.

A remarkable property of an \emph{h-convex } hypersurface $M=\partial\Omega$ in $\mathbb{H}^{n+1}$ is that its inner radius and outer radius are comparable. Recall that the inner radius $\rho_-$ and outer radius $\rho_+$ of a bounded domain $\Omega$ are defined by
\begin{equation*}
  \rho_-= ~\sup\{\rho: B_{\rho}(p)\subset\Omega \text{\rm\ for some }p\in{\mathbb H}^{n+1}\}
  \end{equation*}
  and
  \begin{equation*}
  \rho_+= ~\inf\{\rho: \Omega\subset B_{\rho}(p)\text{\rm\ for some }p\in{\mathbb H}^{n+1}\},
\end{equation*}
where $B_{\rho}(p)$ denotes the geodesic ball of radius $\rho$ about $p$ in $\mathbb{H}^{n+1}$. The following results can be found in \cite{BoM1999,Bor-Mi2002,Cab-Miq2007,Mak2012}.
\begin{thm}
Let $\Omega$ be a compact \emph{h-convex} domain in $\mathbb{H}^{n+1}$ and denote the center of an inball by $o$ and its inner radius by $\rho_-$. Then we have
\begin{itemize}
  \item[(1)]The maximum of the distance $\mathrm{d}_{\mathbb{H}}(o,\cdot)$ between $o$ and the points on $\partial\Omega$ satisfies
  \begin{equation}\label{s2:inner}
  \max_{p\in\partial\Omega} \mathrm{d}_{\mathbb{H}}(o,p)~\leq~ \rho_-+\ln\frac{(1+\sqrt{\tanh\rho_-/2})^2}{1+\tanh\rho_-/2}<~\rho_-+\ln 2.
\end{equation}
Therefore there exists a constant $c>0$ such that the outer radius
\begin{equation}\label{s2:inner2}
  \rho_+\leq c(\rho_-+\rho_-^{1/2}).
\end{equation}
  \item[(2)]For any interior point $p$ of $\Omega$, and any boundary point $q\in\partial\Omega$,
  \begin{equation}\label{s2:nor}
    Dr_p(\nu(q))~\geq~\tanh(\mathrm{d}_{\mathbb{H}}(p,\partial\Omega)),
  \end{equation}
  where $r_p(x)=\mathrm{d}_{\mathbb{H}}(p,x)$.
\end{itemize}
\end{thm}

For smooth \emph{h-convex} domains in $\mathbb{H}^{n+1}$,  inequalities of Alexandrov-Fenchel type for quermassintegrals were proved by Wang-Xia in \cite{WX}. See also \cite{GWW-2014JDG,LWX-2014} for related Alexandrov-Fenchel type inequalities for the curvature integrals \eqref{s2:CurInt}.
\begin{thm}[\cite{WX}]\label{s2:thm1}
For any smooth bounded domain $\Omega$ in $\mathbb{H}^{n+1}$ with \emph{h-convex} boundary $\partial\Omega$, and $0\leq l<k\leq n$, we have
\begin{equation}\label{s2:AF1}
  W_k(\Omega)~\geq~f_k\circ f_l^{-1}(W_l(\Omega))
\end{equation}
with equality if and only if $\Omega$ is a geodesic ball. Here the function $f_k:[0,\infty)\to \mathbb{R}_+$ is increasing and is defined by $f_k(r)=W_k(B_r)$, with $B_r$ a geodesic ball in $\mathbb{H}^{n+1}$. $f_l^{-1}$ is the inverse function of $f_l$.
\end{thm}

\subsection{Inverse concave functions}
For a smooth symmetric function $F(A)=f(\kappa(A))$, where $A=(A_{ij})\in \mathrm{Sym}(n)$ is a symmetric matrix and $\kappa(A)=(\kappa_1,\cdots,\kappa_n)$ gives the eigenvalues of $A$, we denote by $\dot{F}^{ij}$ and $\ddot{F}^{ij,kl}$ the first and second derivatives of $F$ with respect to the components of its argument, so that
\begin{equation*}
  \frac{\partial}{\partial s}F(A+sB)\bigg|_{s=0}=\dot{F}^{ij}(A)B_{ij}
\end{equation*}
and
\begin{equation*}
  \frac{\partial^2}{\partial s^2}F(A+sB)\bigg|_{s=0}=\ddot{F}^{ij,kl}(A)B_{ij}B_{kl}
\end{equation*}
for any two symmetric matrices $A,B$.  We also use the notation
\begin{equation*}
  \dot{f}^i(\kappa)=\frac{\partial f}{\partial \kappa_i}(\kappa),\quad  \ddot{f}^{ij}(\kappa)=\frac{\partial^2 f}{\partial \kappa_i\partial\kappa_j}(\kappa).
\end{equation*}
for the derivatives of $f$ with respect to $\kappa$. At any diagonal $A$ with distinct eigenvalues, the second derivative $\ddot{F}$ of $F$ in direction $B\in \mathrm{Sym}(n)$ is given in terms of $\dot{f}$ and $\ddot{f}$ by  (see \cite{And1994-2,And2007}):
\begin{equation}\label{s2:F-ddt}
  \ddot{F}^{ij,kl}B_{ij}B_{kl}=\sum_{i,k}\ddot{f}^{ik}B_{ii}B_{kk}+2\sum_{i>k}\frac{\dot{f}^i-\dot{f}^k}{\kappa_i-\kappa_k}B_{ik}^2.
\end{equation}
This formula makes sense as a limit in the case of any repeated values of $\kappa_i$. Since $\Psi(A)=F^{\alpha}(A)$, we will use the same notations and formulas for the derivatives of $\Psi$ and $\psi=f^{\alpha}$.

For any positive definite symmetric matrix $A\in \mathrm{Sym}(n)$ with eigenvalues $\kappa(A)\in\Gamma_+$, define $F_*(A)=F(A^{-1})^{-1}$. Then $F_*(A)=f_*(\kappa(A))$, where $f_*$ is defined in \eqref{s3:f-dual}. Since $f$ is defined on the positive definite cone $\Gamma_+$, the following lemma characterizes the inverse concavity of $f$ and $F$.
 \begin{lem}[\cite{And2007,Andrews-McCoy-Zheng}]\label{s3:lem1}
\begin{itemize}
\item[(i)] $f$ is inverse concave if and only if the following matrix
 \begin{equation}\label{s1:matrx}
   \left(\ddot{f}^{kl}+2\frac {\dot{f}^k}{\kappa_k}\delta_{kl}\right)\geq~0.
 \end{equation}
  \item[(ii)] $F_*$ is concave if and only if $f_*$ is concave;
    \item[(ii)]$F_*$ is inverse concave if and only if
\begin{equation}\label{s3:f-invcon}
   \left(\ddot{f}^{kl}+2\frac {\dot{f}^k}{\kappa_k}\delta_{kl}\right)\geq 0,\quad \mathrm{and }\quad \frac{\dot{f}^k-\dot{f}^l}{\kappa_k-\kappa_l}+\frac{\dot{f}^k}{\kappa_l}+\frac{\dot{f}^l}{\kappa_k}\geq~0,~k\neq l.
\end{equation}
\item[(i)] If $f$ is inverse concave, then
  \begin{equation}\label{s3:f-invcon-0}
    \sum_{i=1}^n\dot{f}^i\kappa_i^2~\geq~f^2.
  \end{equation}
\end{itemize}
 \end{lem}

Since the function $f=E_k^{1/k}$ is inverse concave for all $k=1,\cdots,n$ and has dual function
\begin{equation*}
  f_*(z)=~\left(\frac{E_n(z)}{E_{n-k}(z)}\right)^{1/k},\quad z\in\Gamma_+
\end{equation*}
which vanishes on the boundary of $\Gamma_+$, we have that $f=E_k^{1/k}$, $k=1,\cdots,n$, satisfy the Assumption \ref{s1:assum}. We can also easily see that a convex function $f:\Gamma_+\to\mathbb{R}$ satisfies the Assumption \ref{s1:assum}. Firstly, the inequality \eqref{s1:matrx} is obviously true since $f$ is convex and strictly increasing.  Secondly, the convexity of $f$ implies that
\begin{equation*}
  f(x_1,\cdots,x_n)~\geq~\frac 1n\sum_{i=1}^nx_i.
\end{equation*}
Then the dual function $f_*$ satisfies
\begin{equation*}
  f_*(z_1,\cdots,z_n)=f\left(\frac 1{z_1},\cdots,\frac 1{z_n}\right)^{-1}\leq n\left(\frac 1{z_1}+\cdots+\frac 1{z_n}\right)^{-1}=\frac{E_n(z)}{E_{n-1}(z)}
\end{equation*}
Thus $f_*$ approaches zero on the boundary of $\Gamma_+$.  Other important examples of functions satisfying Assumption \ref{s1:assum} are the power means  $H_r=(\frac 1n\sum_ix_i^r)^{1/r}$, which are inverse-concave for $r\geq -1$, concave for $r\geq 1$, and have $f_*$ approaching zero on $\partial\Gamma_+$ for $r\geq 0$, and so satisfy our requirements for all $r\geq 0$.  More examples can be constructed as follows:  If $G_1$ is homogeneous of degree one, increasing in each argument, and inverse-concave, and $G_2$ satisfies Assumption \ref{s1:assum}, then $F=G_1^\sigma G_2^{1-\sigma}$ satisfies Assumption \ref{s1:assum}for any $0<\sigma<1$ (see \cite{And2007,ALM14} for more examples of inverse concave or convex functions).

\subsection{Evolution equations}
Along the flow
\begin{equation*}
  \frac{\partial}{\partial t}X(x,t)=~(\phi(t)-\Psi(\mathcal{W}(x,t)))\nu(x,t)
\end{equation*}
in hyperbolic space $\mathbb{H}^{n+1}$, we have the following evolution equations (see \cite{LWX-2014}) for the induced metric $g_{ij}$, unit outward normal $\nu$, induced area element $d\mu_t$ and Weingarten matrix $\mathcal{W}=(h_i^j)$ of $M_t=X(M^n,t)$:
\begin{align}
  \frac{\partial}{\partial t}g_{ij} =&~ 2(\phi(t)-\Psi)h_{ij} \label{evl-g}\displaybreak[0]\\
  \frac{\partial}{\partial t}\nu=&~\nabla \Psi\label{evl-nu}\displaybreak[0]\\
   \frac{\partial}{\partial t}d\mu_t=& ~nE_1 (\phi(t)-\Psi)d\mu_t\label{evl-dmu}\displaybreak[0]\\
    \frac{\partial}{\partial t}h_i^j=&~\nabla^j\nabla_i\Psi+(\Psi-\phi(t))(h_i^kh_k^j-\delta_i^j)\label{s2:evl-h}
\end{align}
where $\nabla$ denotes the Levi-Civita connection with respect to the induced metric $g_{ij}$ on $M_t$. As an immediate consequence of \eqref{s2:evl-h}, we have that the curvature function $\Psi=\Psi(\mathcal{W})$ evolves by
\begin{equation}\label{s2:evl-F}
  \frac{\partial}{\partial t}\Psi=~\dot{\Psi}^{kl}\nabla^k\nabla_l\Psi+(\Psi-\phi(t))(\dot{\Psi}^{ij}h_i^kh_k^j-\dot{\Psi}^{ij}\delta_i^j),
\end{equation}
where $\dot{\Psi}^{kl}$ denotes the derivatives of $\Psi$ with respect to the components of $\mathcal{W}=(h_i^j)$. Throughout this paper we will always evaluate the derivatives of $\Psi=F^{\alpha}$ at $\mathcal{W}=(h_i^j)$ and the derivatives of $\psi=f^{\alpha}$ at $\kappa(\mathcal{W})=(\kappa_1,\cdots,\kappa_n)$.

The following lemma gives a parabolic type equation of $h_i^j$.
\begin{lem}\label{s2:lem-h}
Along the flow \eqref{flow-VMCF}, the Weingarten matrix $h_i^j$ of $M_t$ evolves by
\begin{align}\label{s3:evl-h}
 \frac{\partial}{\partial t}h_i^j  =& \dot{\Psi}^{kl}\nabla_k\nabla_lh_i^j+\ddot{\Psi}^{kl,pq}\nabla_ih_{kl}\nabla^jh_{pq} +(\dot{\Psi}^{kl}h_k^rh_{rl}+\dot{\Psi}^{kl}g_{kl})h_i^j\nonumber\\
   &\quad -\dot{\Psi}^{kl}h_{kl}(h_i^ph_{p}^j+\delta_i^j)+(\Psi-\phi(t))(h_i^kh_{k}^j-\delta_i^j),
\end{align}
where $\Psi=F^{\alpha}$ and $\dot{\Psi}^{kl}, \ddot{\Psi}^{kl,pq}$ denote the derivatives of $\Psi$ with respect to the components of $\mathcal{W}=(h_i^j)$.
\end{lem}
\proof
Firstly, combining the Gauss and Codazzi equation in hyperbolic space gives the following generalized Simons' identity (see \cite{And-chen2014}):
\begin{equation}\label{s2:Sim}
  \nabla_{(i}\nabla_{j)}h_{kl}=\nabla_{(k}\nabla_{l)}h_{ij}+(h_k^ph_{pl}+g_{kl})h_{ij}-h_{kl}(h_i^ph_{pj}+g_{ij}),
\end{equation}
where the brackets denote symmetrisation. Then
\begin{align}\label{s2:nab-Psi}
  \nabla^j\nabla_i\Psi=&~\dot{\Psi}^{kl}\nabla^j\nabla_ih_{kl}+\ddot{\Psi}^{kl,pq}\nabla_ih_{kl}\nabla^jh_{pq} \nonumber\\
  = &~ \dot{\Psi}^{kl}\nabla_{k}\nabla_{l}h_{i}^j+\ddot{\Psi}^{kl,pq}\nabla_ih_{kl}\nabla^jh_{pq}\nonumber\\
  &\quad +\dot{\Psi}^{kl}(h_k^ph_{pl}+g_{kl})h_{i}^j-\dot{\Psi}^{kl}h_{kl}(h_i^ph_{p}^j+\delta_i^j).
\end{align}
The equation \eqref{s3:evl-h} follows from \eqref{s2:evl-h} and \eqref{s2:nab-Psi} immediately.
\endproof
Using the evolution equations \eqref{evl-dmu} and \eqref{s2:evl-h}, we can also derive the evolution equation for the curvature integral defined in \eqref{s2:CurInt}
\begin{align}\label{s2:evl-CurInt}
 \frac{d}{d t}V_{n-k}(\Omega_t) =&\int_{M_t} \left( \frac{\partial}{\partial t}E_kd\mu_t+E_k \frac{\partial}{\partial t}d\mu_t\right) \nonumber\\
  = &\int_{M_t}\left(\frac{\partial E_k}{\partial h_i^j}\nabla^j\nabla_iF^{\alpha}+(F^{\alpha}-\phi(t))( \frac{\partial E_k}{\partial h_i^j}(h_i^kh_k^j-\delta_i^j)+nE_1E_k)\right)d\mu_t\nonumber\\
  = &\int_{M_t}\left((\phi(t)-F^{\alpha})( (n-k)E_{k+1}+kE_{k-1})\right)d\mu_t,
\end{align}
where we used the facts that $ \nabla^j({\partial E_k}/{\partial h_i^j})=0$ and
\begin{equation*}
\frac{\partial E_k}{\partial h_i^j}h_i^kh_k^j=nE_1E_k-(n-k)E_{k+1},\qquad \frac{\partial E_k}{\partial h_i^j}\delta_i^j=kE_{k-1}.
\end{equation*}
By applying induction argument to \eqref{s2:evl-CurInt} and \eqref{s2:Quer-int1}, we have the following evolution equation for the quermassintegrals of $\Omega_t$ along the flow \eqref{flow-VMCF},
\begin{equation}\label{s2:evl-Wk-1}
  \frac d{dt}W_k(\Omega_t)=~\frac {n+1-k}{n+1}\int_{M_t}E_k\left(\phi(t)-F^{\alpha}\right)d\mu_t,\quad k=0,\cdots,n,
\end{equation}
which was also derived in \cite{WX}. Thus for the function $\phi(t)$ defined in \eqref{s1:phit}, the quermassintegral $W_k(\Omega_t)$ remains constant along the flow \eqref{flow-VMCF}.
\begin{rem}
If $\phi(t)$ is defined as
\begin{equation}\label{s2:phit}
  \phi(t)=~\frac 1{|M_t|}\int_{M_t}F^{\alpha}d\mu_t,
\end{equation}
then the volume $|\Omega_t|$ remains constant. The flow \eqref{flow-VMCF} with $\phi(t)$ given by \eqref{s2:phit} is called the volume preserving curvature flow.
\end{rem}

\section{Preserving of h-convexity}\label{sec:h-conv}
In this section, we will use the tensor maximum principle to prove that the \emph{h-convexity}  is preserved along the flow \eqref{flow-VMCF} if $f$ is inverse concave. For the convenience of readers, we include here the statement of the tensor maximum principle, which was first proved by Hamilton \cite{Ha1982} and was generalized by Andrews \cite{And2007}.
\begin{thm}[\cite{And2007}]\label{s3:tensor-mp}
Let $S_{ij}$ be a smooth time-varying symmetric tensor field on a compact manifold $M$, satisfying
\begin{equation}
\frac{\partial}{\partial t}S_{ij}=a^{kl}\nabla_k\nabla_lS_{ij}+u^k\nabla_kS_{ij}+N_{ij},
\end{equation}
where $a^{kl}$ and $u$ are smooth, $\nabla$ is a (possibly time-dependent) smooth symmetric connection, and $a^{kl}$ is positive definite everywhere. Suppose that
\begin{equation}\label{s3:TM2}
  N_{ij}v^iv^j+\sup_{\Lambda}2a^{kl}\left(2\Lambda_k^p\nabla_lS_{ip}v^i-\Lambda_k^p\Lambda_l^qS_{pq}\right)\geq 0
\end{equation}
whenever $S_{ij}\geq 0$ and $S_{ij}v^j=0$. If $S_{ij}$ is positive definite everywhere on $M$ at $t=0$ and on $\partial M$ for $0\leq t\leq T$, then it is positive on $M\times[0,T]$.
\end{thm}

The main result of this section is the following.
\begin{thm}\label{s3:thm-2}
For any power $\alpha>0$, if the initial hypersurface $M_0$ is h-convex and $f$ is inverse concave, then along the flow \eqref{flow-VMCF} in $\mathbb{H}^{n+1}$ the flow hypersurface $M_t$ is strictly h-convex for $t>0$.
\end{thm}
\proof
Denote
\begin{equation*}
  S_{ij}=h_i^j-\delta_i^j.
\end{equation*}
Then the \emph{h-convexity} is equivalent to $S_{ij}\geq 0$. By \eqref{s3:evl-h}, the tensor $S_{ij}$ evolves by
\begin{align}\label{s3:evl-S}
 \frac{\partial}{\partial t}S_{ij}  =&~ \dot{\Psi}^{kl}\nabla_k\nabla_lS_{ij}+\ddot{\Psi}^{kl,pq}\nabla_ih_{kl}\nabla^jh_{pq} +(\dot{\Psi}^{kl}h_k^rh_{rl}+\dot{\Psi}^{kl}g_{kl})S_{ij}\nonumber\\
   &\quad +(\Psi-\phi(t)-\dot{\Psi}^{kl}h_{kl}) (S_{ik}S_{kj}+2S_{ij})\nonumber\\
   &\quad +\dot{\Psi}^{kl}(h_k^rh_{rl}+g_{kl}-2h_{kl})\delta_i^j
\end{align}
To apply the tensor maximum principle in Theorem \ref{s3:tensor-mp}, we need to show the inequality \eqref{s3:TM2} whenever $S_{ij}\geq 0$ and $S_{ij}v^j=0$ (so that $v$ is a null vector of $S$).  Let $(x_0,t_0)$ be the point where $S_{ij}$ has a null vector $v$. By continuity we can assume that $h_{i}^j$ has all eigenvalues distinct and in increasing order at $(x_0,t_0)$, that is $\kappa_n>\kappa_{n-1}>\cdots>\kappa_1$. The null eigenvector condition $S_{ij}v^j=0$ implies that $v=e_1$ and $S_{11}=\kappa_1-1=0$ at $(x_0,t_0)$.

The lower order terms in \eqref{s3:evl-S} involving $S_{ij}$ and $S_{ik}S_{kj}$ satisfy the null vector condition and can be ignored:  In particular for the last term in \eqref{s3:evl-S}, we have
\begin{align}
  \dot{\Psi}^{kl}(h_k^rh_{rl}+g_{kl}-2h_{kl})= & \sum_k\dot{\psi}^k\left(\kappa_k^2+1-2\kappa_k\right) \nonumber\displaybreak[0]\\
=&~  \sum_k\dot{\psi}^k(\kappa_k-1)^2~\geq~0.
\end{align}
Thus it remains to show that
\begin{align}\label{s3:Q1}
  Q_1: =& \ddot{\Psi}^{kl,pq}\nabla_1h_{kl}\nabla_1h_{pq} +2\sup_{\Lambda}\dot{\Psi}^{kl}\left(2\Lambda_k^p\nabla_lS_{1p}-\Lambda_k^p\Lambda_l^qS_{pq}\right)~\geq~0.
\end{align}
Note that $S_{11}=0$ and $\nabla_kS_{11}=0$ at $(x_0,t_0)$, the supremum over $\Lambda$ can be computed exactly as follows:
\begin{align*}
  2\dot{\Psi}^{kl}&\left(2\Lambda_k^p\nabla_lS_{1p}-\Lambda_k^p\Lambda_l^qS_{pq}\right)  \\
  &\quad =2\sum_{k=1}^n\sum_{p=2}^n\dot{\psi}^k\left(2\Lambda_k^p\nabla_kS_{1p}-(\Lambda_k^p)^2S_{pp}\right)\\
  &\quad =2 \sum_{k=1}^n\sum_{p=2}^n\dot{\psi}^k\left(\frac{(\nabla_kS_{1p})^2}{S_{pp}}-\left(\Lambda_k^p-\frac{\nabla_kS_{1p}}{S_{pp}}\right)^2S_{pp}\right).
\end{align*}
It follows that the supremum is obtained by choosing $\Lambda_k^p=\frac{\nabla_kS_{1p}}{S_{pp}}$. The required inequality for $Q_1$ becomes:
\begin{align}\label{s3:Q1-1}
  Q_1=& \ddot{\Psi}^{kl,pq}\nabla_1h_{kl}\nabla_1h_{pq} +2 \sum_{k=1}^n\sum_{p=2}^n\dot{\psi}^k\frac{(\nabla_kS_{1p})^2}{S_{pp}}~\geq~0.
\end{align}
Using \eqref{s2:F-ddt} to express the second derivatives of $\Psi$ and noting that $\psi=f^{\alpha}$, $\nabla_1S_{1p}=\nabla_1h_{1p}=\nabla_ph_{11}=0$ at $(x_0,t_0)$, we have
\begin{align}\label{s3:Q1-0}
  Q_1=& \alpha f^{\alpha-1}\ddot{f}^{kl}\nabla_1h_{kk}\nabla_1h_{ll} +\alpha(\alpha-1)f^{\alpha-2}(\nabla_1F)^2\nonumber\\
  &\quad +2\alpha f^{\alpha-1}\sum_{k>l}\frac{\dot{f}^k-\dot{f}^l}{\kappa_k-\kappa_l}(\nabla_1h_{kl})^2+2 \alpha f^{\alpha-1}\sum_{k>1,l>1}\frac{\dot{f}^k}{\kappa_{l}-1}(\nabla_1h_{kl})^2.
\end{align}
To make use the inverse concavity of $f$, let $\tau_i=1/{\kappa_i}$ and $f_*(\tau)=f(\kappa)^{-1}$. We compute that
 \begin{align*}
  \dot{f}^k =& f_*^{-2}\frac{\partial f_*}{\partial \tau_k}\frac 1{\kappa_k^2} \\
  \ddot{f}^{kl} =& -f_*^{-2}\frac{\partial^2 f_*}{\partial \tau_k\partial \tau_l}\frac 1{\kappa_k^2\kappa_l^2}+2f_*^{-3}\frac{\partial f_*}{\partial \tau_k}\frac 1{\kappa_k^2}\frac{\partial f_*}{\partial \tau_l}\frac 1{\kappa_l^2}-2f_*^{-2}\frac{\partial f_*}{\partial \tau_k}\frac 1{\kappa_k^3}\delta_{kl}\\
  =&-f_*^{-2}\frac{\partial^2 f_*}{\partial \tau_k\partial \tau_l}\frac 1{\kappa_k^2\kappa_l^2}+2f^{-1}\dot{f}^k\dot{f}^l-2\frac{\dot{f}^k}{\kappa_k}\delta_{kl}.
\end{align*}
By the concavity of $f_*$, the first term of \eqref{s3:Q1-0} can be estimated as
\begin{align*}
  \alpha f^{\alpha-1}\ddot{f}^{kl,pq}\nabla_1h_{kl}\nabla_1h_{pq} \geq &  2\alpha f^{\alpha-1}\left(f^{-1}(\nabla_1F)^2-\sum_k\frac{\dot{f}^k}{\kappa_k}(\nabla_1h_{kk})^2\right)
\end{align*}
Then
\begin{align*}
  \frac{Q_1}{ \alpha f^{\alpha-1}}\geq &(\alpha+1)f^{-1}(\nabla_1F)^2-2\sum_k\frac{\dot{f}^k}{\kappa_k}(\nabla_1h_{kk})^2\\
  &\quad +2\sum_{k>l}\frac{\dot{f}^k-\dot{f}^l}{\kappa_k-\kappa_l}(\nabla_1h_{kl})^2+2\sum_{k>1,l>1}\frac{\dot{f}^k}{\kappa_{l}-1}(\nabla_1h_{kl})^2\displaybreak[0]\\
  \geq &(\alpha+1)f^{-1}(\nabla_1F)^2-2\sum_{k>1}\frac{\dot{f}^k}{\kappa_k}(\nabla_1h_{kk})^2\displaybreak[0]\\
  &\quad -2\sum_{k\neq l>1}\frac{\dot{f}^k}{\kappa_l}(\nabla_1h_{kl})^2+2\sum_{k>1,l>1}\frac{\dot{f}^k}{\kappa_{l}-1}(\nabla_1h_{kl})^2\displaybreak[0]\\
 =&(\alpha+1)f^{-1}(\nabla_1F)^2+2\sum_{k>1,l>1}\left(\frac{\dot{f}^k}{\kappa_{l}-1}-\frac{\dot{f}^k}{\kappa_l}\right)(\nabla_1h_{kl})^2\\
 \geq&~0
  \end{align*}
for all $\alpha>0$, where we used the second inequality of \eqref{s3:f-invcon} and the fact $\nabla_kh_{11}=0$. The tensor maximum principle implies that the \emph{h-convexity} is preserved along the flow \eqref{flow-VMCF}.

Finally we show that $M_t$ is strictly \emph{h-convex} for $t>0$. If this is not true,  there exists some interior point $(x_0,t_0)$ such that the smallest principal curvature is $1$.  By the strong maximum principle there exists a parallel vector field $v$ such that $S_{ij}v^iv^j=0$ on $M_{t_0}$. Then the smallest principal curvature is 1 on $M_{t_0}$ everywhere. On the other hand, a standard argument shows that on any closed hypersurface in $\mathbb{H}^{n+1}$, there exists at least one point where all the principal curvatures are strictly bigger than one. This contradiction completes the proof of Theorem \ref{s3:thm-2}.
\endproof

\begin{rem}\label{s3:rem}
The above argument implies that \emph{h-convexity} is also preserved along the flow
\begin{equation*}
  \frac{\partial}{\partial t}X(x,t)=~(\phi(t)-\Psi(\mathcal{W}(x,t)))\nu(x,t)
\end{equation*}
in hyperbolic space $\mathbb{H}^{n+1}$ with $\Psi=\psi(f)$, where $f$ is an inverse concave function and $\psi:[0,+\infty)\to \mathbb{R}_+$ satisfies $\psi'(r)>0$ and $\psi''(r)r+2\psi'(r)\geq 0$ for all $r>0$.
\end{rem}

\section{Upper bound of $F$}\label{sec:ub-F}
In this section, we will prove that $F$ is uniformly bounded from above along the flow \eqref{flow-VMCF}. Firstly, the preservation of $W_k(\Omega_t)$ and the\emph{ h-convexity} of $M_t=\partial \Omega_t$ imply uniform two-sided bounds on inner radius and outer radius of $\Omega_t$.
\begin{lem}
Denote by $\rho_-(t)$ and $\rho_+(t)$ the inner and outer radii of the domain $\Omega_t$ enclosed by $M_t$. Then there exist positive constants $c_1,c_2$ depending only on $n,k,M_0$ such that
\begin{equation}\label{s4:io-radius1}
  0<c_1\leq \rho_-(t)\leq \rho_+(t)\leq c_2
\end{equation}
for all time $t\in [0,T)$.
\end{lem}
\proof
On the one hand, since $W_k(\Omega_t)=W_k(\Omega_0)$, we have
\begin{equation*}
  W_k(B_{\rho_+(t)})\geq ~W_k(\Omega_t)=W_k(\Omega_0),
\end{equation*}
where $B_{\rho_+(t)}$ is the geodesic ball of radius $\rho_+(t)$ that encloses $\Omega_t$. Thus
\begin{equation*}
  \rho_+(t)\geq f_k^{-1}(W_k(\Omega_0))>0,
\end{equation*}
where $f_k^{-1}$ is the inverse function of $f_k(r)=W_k(B_r)$. Similarly, $\rho_-(t)\leq f_k^{-1}(W_k(\Omega_0))$. Since each $M_t$ is \emph{h-convex}, the estimate \eqref{s4:io-radius1} follows by the inequality \eqref{s2:inner2}.
\endproof

By \eqref{s4:io-radius1}, the inner radius of $\Omega_t$ is bounded below by a positive constant $c_1$. This implies that there exists a geodesic ball of radius $c_1$ contained in $\Omega_t$ for each $t\in [0,T)$. The following lemma shows the existence of a geodesic ball with fixed center enclosed by the flow hypersurfaces on a suitable time interval.

\begin{lem}\label{s5:lem-inball}
Let $M_t$ be a smooth h-convex solution of \eqref{flow-VMCF} on $[0,T)$ with global term $\phi(t)$ given by \eqref{s1:phit}. For any $t_0\in [0,T)$, let $B(p_0,\rho_0)$ be the inball of $\Omega_{t_0}$, where $\rho_0=\rho_-(t_0)$. Then
\begin{equation}\label{s5:inball-eqn1}
  B(p_0,\rho_0/2)\subset \Omega_t,\quad t\in [t_0, \min\{T,t_0+\tau\})
\end{equation}
for some $\tau$ depending only on $n,\alpha,k,\Omega_0$.
\end{lem}
\proof
Given $p_0$, we denote by $r_{p_0}$ the distance function to $p_0$ in $\mathbb{H}^{n+1}$ and by $\partial_r=\partial_{r_{p_0}}$ the gradient vector of $r_{p_0}$. For any $x\in M_t$,
\begin{align}\label{s5:inball-1}
  \frac{\partial}{\partial t} \sinh^2r_{p_0}(x)=& 2\langle \sinh r_{p_0}(x)\partial r, \frac{\partial}{\partial t}(\sinh r_{p_0}(x)\partial_r) \rangle\nonumber\\
  = &2 \sinh r_{p_0}(x)\cosh r_{p_0}(x)(\phi(t)-F^{\alpha}(x,t))\langle \partial_r,\nu\rangle,
\end{align}
where we used the conformal property of the vector field $\sinh r\partial r$, i.e.,
 \begin{equation}\label{s5:conf}
   \langle \bar{\nabla}_X(\sinh r\partial_r), Y\rangle=\cosh r \langle  X,Y\rangle
 \end{equation}
for any tangential vector fields $X,Y$ in $\mathbb{H}^{n+1}$ (see, e.g.,\cite{GL2015}). It follows from \eqref{s5:inball-1} that
 \begin{align*}
  \frac{\partial}{\partial t} r_{p_0}(x)=& (\phi(t)-F^{\alpha}(x,t))\langle \partial r,\nu\rangle\geq~-F^{\alpha}(x,t)\langle \partial_r,\nu\rangle,
\end{align*}
 since $\phi(t)>0$ and $\langle \partial_r,\nu\rangle>0$ on $M_t$. Denote $r(t)=\min_{M_t}r_{p_0}(x)$. At the minimum point, we have $\langle \partial_r,\nu\rangle=1$ and $\kappa_i\leq \coth r(t)$. Then $F\leq \coth r(t)$ at the minimum point and
\begin{equation}\label{s5:inball-3}
  \frac d{dt}r(t)\geq -\coth\left(r(t)\right)^{\alpha}.
\end{equation}
Note that $0<c_1\leq r(t)\leq 2\rho_+(t)\leq 2c_2$, where $c_1,c_2$ are the constants in \eqref{s4:io-radius1}. Then
\begin{equation}\label{s5:c3}
  \coth^{\alpha-1}r(t)~\leq ~\max\{\coth^{\alpha-1}(2c_2),\coth^{\alpha-1}(c_1)\}~=:~c_3.
\end{equation}
We deduce from \eqref{s5:inball-3} and \eqref{s5:c3} that
\begin{equation*}
  \frac d{dt}r(t)\geq -c_3\coth r(t),
\end{equation*}
from which we solve that
\begin{equation*}
  \cosh r(t)~\geq~\cosh r(0)\exp \{-c_3t\}.
\end{equation*}
In particular,
\begin{equation*}
  r(t)\geq~\frac{r(0)}2=~\frac{\rho_0}2
\end{equation*}
provided that
\begin{equation*}
  t-t_0\leq~\frac 1{c_3}\ln\frac{\cosh r(0)}{\cosh \frac{r(0)}2}~=:~\tau
\end{equation*}
 which depends only on $n,\alpha, k, \Omega_0$. Then $B(p_0,\rho_0/2)\subset \Omega_t$ for $t\in [t_0, \min\{T,t_0+\tau\})$.
\endproof

Consider the support function $u(x,t)=\sinh r_{p_0}(x)\langle \partial_{r_{p_0}},\nu\rangle $ of $M_t$ with respect to the point $p_0$. Then by \eqref{s2:nor} and \eqref{s5:inball-eqn1},
\begin{equation}\label{s5:sup-1}
  u(x,t)~\geq ~\sinh(\frac{\rho_0}2)\tanh(\frac{\rho_0}2)~=:~2c
\end{equation}
on $M_t$ for any $t\in[t_0,\min\{T,t_0+\tau\})$. On the other hand, the estimate \eqref{s4:io-radius1} implies that $u(x,t)\leq \sinh(2c_2)$ on $M_t$ for all $t\in[0,T)$.
\begin{lem}
The support function $u(x,t)$ evolves by
\begin{align}\label{s5:evl-u}
   \frac{\partial}{\partial t}u =& \dot{\Psi}^{kl}\nabla_k\nabla_lu+\cosh r_{p_0}(x)\left(\phi(t)-\Psi-\dot{\Psi}^{kl}h_{kl}\right)+\dot{\Psi}^{ij}h_i^kh_{kj}u.
\end{align}
\end{lem}
 \proof
 Firstly, by \eqref{s5:conf} and \eqref{evl-nu},
 \begin{equation}\label{s5:evl-u1}
   \frac{\partial}{\partial t}u=\cosh r_{p_0}(x) \left(\phi(t)-\Psi\right)+\sinh r_{p_0}(x)\langle\partial r,\nabla \Psi\rangle.
 \end{equation}
 Secondly, the spatial derivatives of $u$ can also be computed using \eqref{s5:conf}:
 \begin{align}
   \nabla_ju=&~\sinh r_{p_0}(x)\langle\partial r,h_i^k\partial_k\rangle \\
   \nabla_i\nabla_ju =&~\cosh r_{p_0}(x)h_{ij}+ \sinh r_{p_0}(x)\langle\partial_r, \nabla^kh_{ij}\partial_k-h_i^kh_{kj}\nu\rangle.\label{s5:evl-u3}
 \end{align}
 Then the evolution equation \eqref{s5:evl-u} follows by a direct computation using \eqref{s5:evl-u1} --\eqref{s5:evl-u3}.
 \endproof
Now we can use the technique that was first introduced by Tso \cite{Tso85} to prove the upper bound of $F$ along the flow \eqref{flow-VMCF}.

\begin{thm}\label{Ek-ub}
Let $M_t$ be a smooth h-convex solution of \eqref{flow-VMCF} on $[0,T)$ with global term $\phi(t)$ given by \eqref{s1:phit}. Then we have $\max_{M_t}F\leq C$ for any $t\in [0,T)$, where $C$ depends on $n,k,\alpha,M_0$ but not on $T$.
\end{thm}
\proof
For any given $t_0\in [0,T)$,  define the auxiliary function
\begin{equation*}
  W(x,t)=\frac {\Psi(x,t)}{u(x,t)-c}
\end{equation*}
which is well-defined for all $t\in [t_0,\min\{T,t_0+\tau\})$  by \eqref{s5:sup-1}. Combining \eqref{s2:evl-F} and \eqref{s5:evl-u}, we can compute the evolution equation of the function $W$
\begin{align}\label{s5:evl-W-1}
   \frac{\partial}{\partial t}W= &\dot{\Psi}^{ij}\left(\nabla_j\nabla_iW+\frac 2{u-c}\nabla_iu\nabla_jW\right)\nonumber \\
  &\quad-\frac{\phi(t)}{u-c}\left( \dot{\Psi}^{ij}(h_i^kh_{k}^j-\delta_i^j)+W\cosh r_{p_0}(x)\right)\nonumber\\
  &\quad +\frac{\Psi}{(u-c)^2}(\Psi+\dot{\Psi}^{kl}h_{kl})\cosh r_{p_0}(x)-\frac{c\Psi}{(u-c)^2}\dot{\Psi}^{ij}h_i^kh_{k}^j-W\dot{\Psi}^{ij}\delta_i^j.
\end{align}
The \emph{h-convexity} of $M_t$ implies that $\dot{\Psi}^{ij}(h_i^kh_{k}^j-\delta_i^j)\geq 0$. So the terms involving $\phi(t)$ in \eqref{s5:evl-W-1} are non-positive. Since $\Psi=F^{\alpha}$ where $F$ is inverse concave, we have $\Psi+\dot{\Psi}^{kl}h_{kl}=(1+\alpha)\Psi$ and $\dot{\Psi}^{ij}h_i^kh_{k}^j\geq \alpha F^{\alpha+1}$. The last term of \eqref{s5:evl-W-1} is obviously non-positive. Therefore,
\begin{align}\label{s5:evl-W-1-2}
   \frac{\partial}{\partial t}W\leq & ~\dot{\Psi}^{ij}\left(\nabla_j\nabla_iW+\frac 2{u-c}\nabla_iu\nabla_jW\right)\nonumber \\
  &\quad +(\alpha+1)W^2\cosh r_{p_0}(x)-\alpha cW^2F.
\end{align}
Using \eqref{s5:sup-1} and the upper bound $r_{p_0}(x)\leq 2c_2$, we obtain the following estimate
\begin{align}\label{evl-W-2}
   \frac{\partial}{\partial t}W\leq &~ \dot{\Psi}^{ij}\left(\nabla_j\nabla_iW+\frac 2{u-c}\nabla_iu\nabla_jW\right)\nonumber \\
  &\quad +W^2\left((\alpha+1)\cosh (2c_2)-\alpha c^{1+\frac 1{\alpha}}W^{1/{\alpha}}\right)
\end{align}
holds on $[t_0,\min\{T,t_0+\tau\})$. Let $\tilde{W}(t)=\sup_{M_t}W(\cdot,t)$. Then \eqref{evl-W-2} implies that
\begin{equation*}
  \frac d{dt}\tilde{W}(t)\leq \tilde{W}^2\left((\alpha+1)\cosh (2c_2)-\alpha c^{1+\frac 1{\alpha}}\tilde{W}^{1/{\alpha}}\right)
\end{equation*}
from which it follows using the maximum principle that
\begin{equation}\label{evl-W-4}
  \tilde{W}(t)\leq \max\left\{ \left(\frac {2(1+\alpha)\cosh(2c_2)}{\alpha}\right)^{\alpha}c^{-(\alpha+1)}, \left(\frac {2}{1+\alpha}\right)^{\frac{\alpha}{1+\alpha}}c^{-1}(t-t_0)^{-\frac{\alpha}{1+\alpha}}\right\}.
\end{equation}
Then the upper bound on $F$ follows from \eqref{evl-W-4} and the facts that
\begin{equation*}
  c=\frac 12\sinh(\frac{\rho_0}2)\tanh(\frac{\rho_0}2)\geq \frac 12\sinh(\frac{c_1}2)\tanh(\frac{c_1}2)
\end{equation*}
and $u-c\leq 2c_2$, where $c_1,c_2$ are constants in \eqref{s4:io-radius1} depending only on $n,k,M_0$.
\endproof

As an corollary of the upper bound of $F$ and the \emph{h-convexity} of $M_t$, there exist constants $c_4,c_5$ depend only on $n,k,\alpha,M_0$ such that
\begin{equation}\label{s5:phit-bd}
  c_4\leq~\phi(t)\leq c_5
\end{equation}
on $[0,T)$.

\section{Long-time existence}\label{sec:LTE}

Until now, we only used the fact that $f$ is inverse concave. The upper bound on $F$ proved in \S \ref{sec:ub-F} implies that the dual function $f_*$ of $f$ is bounded from below by a positive constant. Applying condition (vi) in Assumption \ref{s1:assum} that $f_*$ approaches zero on the boundary of $\Gamma_+$, there exists a positive constant $C$ such that $1/{\kappa_i}\geq C$ for all $i=1,\cdots,n$, which implies a uniform upper bound on the Weingarten matrix $\mathcal{W}=(h_{i}^j)$ along the flow \eqref{flow-VMCF} for all $t\in[0,T)$.
\begin{lem}\label{s5:lem-kapp}
There exists a constant $C>0$ depending only on $n,k,,\alpha,M_0$ such that along the flow \eqref{flow-VMCF}, the principal curvatures $\kappa=(\kappa_1,\cdots,\kappa_i)$ of the solution $M_t$ satisfy
\begin{equation*}
  1\leq \kappa_i\leq C,\quad i=1,\cdots,n
\end{equation*}
for all $t\in[0,T)$.
\end{lem}
In the following, we will derive higher order estimates on the solution $M_t$ of the flow \eqref{flow-VMCF} and prove that the solution $M_t$ exists for all time $t\in[0,\infty)$.

Let us denote by $\mathbb{R}^{1,n+1}$ the Minkowski spacetime, that is the vector space $\mathbb{R}^{n+2}$ endowed with the Minkowski spacetime metric $\langle \cdot,\cdot\rangle$ given by
\begin{equation*}
  \langle X,X\rangle~=~-X_0^2+\sum_{i=1}^nX_i^2
\end{equation*}
for any vector  $X=(X_0,X_1,\cdots,X_n)\in \mathbb{R}^{n+2}$. The hyperbolic space $\mathbb{H}^{n+1}$ is then
\begin{equation*}
  \mathbb{H}^{n+1}=~\{X\in \mathbb{R}^{1,n+1},~~\langle X,X\rangle=-1,~X_0>0\}
\end{equation*}
An embedding $X:M^n\to \mathbb{H}^{n+1}$ induces an embedding $Y:M^n\to B_1(0)\subset \mathbb{R}^{n+1}$ by
\begin{equation}\label{s5:XY}
  X~=~\frac{(1,Y)}{\sqrt{1-|Y|^2}},
\end{equation}
where $Y\in B_1(0)\subset \mathbb{R}^{n+1}$.  Let $\{x_i\}, i=1,\cdots,n$ be a local coordinate system on $M$ and $\{\partial_i\}$ be the corresponding coordinate vectors. Then
\begin{equation}\label{s5:diX}
  \partial_iX=~X_*(\partial_i)=~\frac{(0,\partial_iY)}{\sqrt{1-|Y|^2}}+\frac{X}{1-|Y|^2}\langle Y,\partial_iY\rangle
\end{equation}
and
\begin{align*}
  \partial_i\partial_jX= &~ \frac{(0,\partial_i\partial_jY)}{\sqrt{1-|Y|^2}}+\frac{(0,\partial_jY)}{(1-|Y|^2)^{3/2}}\langle Y,\partial_iY\rangle+\partial_iX\frac{\langle Y,\partial_iY\rangle} {1-|Y|^2}\nonumber\\
  & \quad+X\partial_i\left(\frac{\langle Y,\partial_iY\rangle} {1-|Y|^2}\right).
\end{align*}
Let $\nu\in T\mathbb{H}^{n+1}, h_{ij}^X$ and $N\in \mathbb{R}^{n+1}, h_{ij}^Y$ be the unit normal vectors and the second fundamental forms of $X(M^n)\subset \mathbb{H}^{n+1}$ and $Y(M^n)\subset \mathbb{R}^{n+1}$ respectively. Note that $\langle \nu,X\rangle=\langle \nu,\partial_iX\rangle=0$.  Taking the inner product of \eqref{s5:diX} with $\nu$, we also have
\begin{equation}\label{s5:nu-1}
  \langle \nu, (0,\partial_iY)\rangle=0
\end{equation}
Therefore,
\begin{align}\label{s5:h-1}
  h_{ij}^X= & -\langle \partial_i\partial_jX,\nu\rangle= -\frac{1}{\sqrt{1-|Y|^2}}\langle (0,\partial_i\partial_jY),\nu\rangle\nonumber\\
   =& ~-\frac{\langle \partial_i\partial_jY,N\rangle}{\sqrt{1-|Y|^2}} \langle (0,N),\nu\rangle =\frac{h_{ij}^Y}{\sqrt{1-|Y|^2}}\langle (0,N),\nu\rangle
\end{align}
By \eqref{s5:nu-1}, we can write $\nu=a_1(0,N)+a_2(1,0)$ where $a_1,a_2$ can be computed as follows:
\begin{align*}
  1= &\langle \nu,\nu\rangle=~a_1^2-a_2^2 \\
  0= &\langle \nu,X\rangle=~\frac 1{\sqrt{1-|Y|^2}}(a_1\langle N,Y\rangle-a_2)
\end{align*}
from which we deduce that
\begin{equation*}
\langle (0,N),\nu\rangle =~a_1=\frac 1{\sqrt{1-\langle N,Y\rangle^2}}.
\end{equation*}
Substituting this into \eqref{s5:h-1} yields
\begin{align}\label{s5:h-2}
  h_{ij}^X= & \frac{h_{ij}^Y}{\sqrt{(1-|Y|^2)(1-\langle N,Y\rangle^2)}}.
\end{align}
From \eqref{s5:diX}, we also have that the induced metrics $g_{ij}^X$ and $g_{ij}^Y$ of $X(M^n)\subset \mathbb{H}^{n+1}$ and $Y(M^n)\subset \mathbb{R}^{n+1}$ are related by
\begin{align*}
  g_{ij}^X =& \langle \partial_iX,\partial_jX\rangle =\frac{\langle \partial_iY,\partial_jY\rangle }{1-|Y|^2}+\frac{\langle Y,\partial_iY\rangle\langle  Y,\partial_jY\rangle }{(1-|Y|^2)^2} \nonumber\\
   =&\frac 1 {1-|Y|^2}\left(g_{ij}^Y+\frac{\langle Y,\partial_iY\rangle\langle  Y,\partial_jY\rangle }{(1-|Y|^2)}\right)
\end{align*}

Suppose $X:M^n\times[0,T)\to \mathbb{H}^{n+1}$ is a solution to the flow \eqref{flow-VMCF}. We next derive the evolution equation of the corresponding $Y: M^n\times[0,T)\to \mathbb{R}^{n+1}$ related by \eqref{s5:XY}. Denote by $\mathcal{W}^X=(h_{ik}^Xg_X^{kj})$ the Weingarten matrix of $X(M^n,t)\subset \mathbb{H}^{n+1}$, where $(g_X^{kj})$ denotes the inverse matrix of the induce metric $g_{kj}^X$. Then
\begin{align*}
  \phi(t)-F^{\alpha}(\mathcal{W}^X) =& ~\langle \partial_tX,\nu\rangle \\
  = & ~\frac{\langle(0,\partial_tY),\nu\rangle}{\sqrt{1-|Y|^2}}+\langle X,\nu\rangle \frac {\langle Y,\partial_tY\rangle}{1-|Y|^2}\\
  =&~\langle \partial_tY,N\rangle \frac{\langle (0,N),\nu\rangle}{\sqrt{1-|Y|^2}}\\
  =&~\langle \partial_tY,N\rangle \frac 1{\sqrt{(1-|Y|^2)(1-\langle N,Y\rangle^2)}}
\end{align*}
where we used the facts $\langle X,\nu\rangle=0$ and $\langle (0,\partial_iY),\nu\rangle=0$ in the third equality. Thus up to a tangential diffeomorphism, $Y: M^n\times[0,T)\to \mathbb{R}^{n+1}$ satisfies the following evolution equation:
\begin{equation}\label{s5:Y-evl}
 \partial_tY=~\sqrt{(1-|Y|^2)(1-\langle N,Y\rangle^2)}\left(\phi(t)-F^{\alpha}(\mathcal{W}^X) \right)N.
\end{equation}
Here $\mathcal{W}^X$ is the Weingarten matrix of $X(M^n,t)\subset \mathbb{H}^{n+1}$, which we next relate to the geometry of $Y$:  In local coordinates, the inverse matrix $\mathcal{W}_X^{-1}$ of $\mathcal{W}^X$ satisfies
\begin{align}\label{s5:W-inv}
  (\mathcal{W}_X^{-1})_{ij} =& (h^{-1}_X)^{jk}g_{ki}^X \nonumber\\
   =&(h^{-1}_Y)^{kj}\left(g_{ki}^Y+\frac{\langle  Y,\partial_iY\rangle\langle Y,\partial_kY\rangle }{(1-|Y|^2)}\right)\sqrt{\frac{1-\langle N,Y\rangle^2}{1-|Y|^2}}
\end{align}
By the estimate \eqref{s4:io-radius1}, $X$ stays in a bounded subset of $\mathbb{H}^{n+1}$. Then there exists a positive constant $c>0$ depending only on $n,k,M_0$ such that
\begin{equation*}
  0<c\leq 1-|Y|^2\leq 1,\quad 0<c\leq 1-\langle N,Y\rangle^2 \leq 1.
\end{equation*}
Since each $M_t$ is \emph{h-convex } in $\mathbb{H}^{n+1}$, the equation \eqref{s5:h-2} implies that each $Y_t=Y(M^n,t)$ is strictly convex in $\mathbb{R}^{n+1}$.  We can parametrise $Y_t$ using the Gauss map and the support function $s(z):=\langle Y(N^{-1}(z)),z\rangle$, where $N^{-1}:\mathbb{S}^n\to M^n$ is the inverse map of the Gauss map which exists due to the strict convexity of $Y_t$. Then $Y$ is given by the embedding $Y:\mathbb{S}^n\to \mathbb{R}^{n+1}$ with (cf. \cite{And2000,Andrews-McCoy-Zheng})
\begin{equation*}
  Y(z)=~s(z)z+\bar{\nabla} s
\end{equation*}
where $\bar{\nabla}$ is the gradient with respect to the standard round metric $\bar{g}$  on $\mathbb{S}^n$. The derivative of this map is given by
\begin{equation*}
  \partial_iY=~\tau_{ik}\bar{g}^{kl}\partial_lz
\end{equation*}
in local coordinates, where $\tau_{ij}$ is given as follows
\begin{equation*}
  \tau_{ij}=~\bar{\nabla}_i\bar{\nabla}_js+s\bar{g}_{ij}.
\end{equation*}
The eigenvalues $\tau_i$ of the matrix $\tau_{ij}$ with respect to the metric $\bar{g}$ are the inverse of the principal curvatures $\kappa_i$, i.e., $\tau_i=1/{\kappa_i}$, and are called the principal radii of curvature. Then
\begin{equation*}
  g^Y_{ki}=~\langle \partial_kY,\partial_iY\rangle=\tau_{kp}\bar{g}^{pq}\tau_{qi}
\end{equation*}
and
\begin{equation*}
  \frac{\langle Y,\partial_iY\rangle\langle Y,\partial_kY\rangle }{(1-|Y|^2)}=~\tau_{kp}\frac{\langle \bar{g}^{pa}\bar{\nabla}_as,\bar{g}^{qb}\bar{\nabla}_bs\rangle}{1-s^2-|\bar{\nabla}s|^2}\tau_{qi}
\end{equation*}
We can express \eqref{s5:W-inv} using $s,\bar{\nabla}s$ and the matrix $\tau_{ij}$ as follows:
\begin{align}\label{s5:W-inv-2}
  (\mathcal{W}_X^{-1})_{ij}    =&(h^{-1}_Y)^{kj}\left(g_{ki}^Y+\frac{\langle  Y,\partial_iY\rangle\langle  Y,\partial_kY\rangle }{(1-|Y|^2)}\right)\sqrt{\frac{1-\langle N,Y\rangle^2}{1-|Y|^2}}\nonumber\\
  =&\tau_{lj}(\tau^{-1})^{ls}\bar{g}_{st}(\tau^{-1})^{tk}\tau_{kp}\left(\bar{g}^{pq}+\frac{\langle \bar{g}^{pa}\bar{\nabla}_as,\bar{g}^{qb}\bar{\nabla}_bs\rangle}{1-s^2-|\bar{\nabla}s|^2}\right)\tau_{qi}\sqrt{\frac{1-s^2}{1-s^2-|\bar{\nabla}s|^2}}\nonumber\\
  =&\left(\bar{g}^{jq}+\frac{\langle \bar{g}^{ja}\bar{\nabla}_as,\bar{g}^{qb}\bar{\nabla}_bs\rangle}{1-s^2-|\bar{\nabla}s|^2}\right)\tau_{qi}\sqrt{\frac{1-s^2}{1-s^2-|\bar{\nabla}s|^2}}
\end{align}
where $\tau^{-1}$ denotes the inverse matrix of $\tau_{ij}$. The solution of the evolution equation \eqref{s5:Y-evl} is then given up to a tangential diffeomorphism by solving the following scalar parabolic equation
\begin{equation}\label{s5:Y-evl-Gau}
  \partial_ts=~\sqrt{(1-s^2-|\bar{\nabla}s|^2)(1-s^2)}\left(\phi(t)-F_{*}^{-\alpha}((\mathcal{W}_X^{-1})_{ij} )\right)
\end{equation}
for the support function $s(z,t)$, where $(\mathcal{W}_X^{-1})_{ij} $ is the matrix given in \eqref{s5:W-inv-2} in terms of $s,\bar{\nabla}s$ and the matrix $\tau_{ij}=\bar{\nabla}_i\bar{\nabla}_js+s\bar{g}_{ij}$.

By \eqref{s4:io-radius1} and Lemma \ref{s5:lem-kapp}, we already have uniform $C^2$ estimates on the support function $s(z,t)$. Denote the right hand side of \eqref{s5:Y-evl-Gau} by $G(\bar{\nabla}^2s, \bar{\nabla}s, s,z,t)$. Then
\begin{align*}
  \dot{G}^{ij}=&~\frac{\partial G}{\partial (\bar{\nabla}^2_{ij}s)}=~\alpha F_*^{-\alpha-1}\sqrt{(1-s^2-|\bar{\nabla}s|^2)(1-s^2)}\dot{F}_*^{pq}\frac{\partial  (\mathcal{W}_X^{-1})_{pq}}{\partial \bar{\nabla}^2_{ij}s}\\
  =&~\alpha F_*^{-\alpha-1}(1-s^2)\dot{F}_*^{pq}\left(\bar{g}^{pi}+\frac{\langle \bar{g}^{pa}\bar{\nabla}_as,\bar{g}^{ib}\bar{\nabla}_bs\rangle}{1-s^2-|\bar{\nabla}s|^2}\right)\bar{g}_{qj},
\end{align*}
and
\begin{align}\label{s6:dd-G}
  \frac 1{\sqrt{(1-s^2-|\bar{\nabla}s|^2)(1-s^2)}}\dot{G}^{ij,kl}=&~\alpha F_*^{-\alpha-1}\dot{F}_*^{pq,rt}\frac{\partial  (\mathcal{W}_X^{-1})_{pq}}{\partial \bar{\nabla}^2_{ij}s}\frac{\partial  (\mathcal{W}_X^{-1})_{rt}}{\partial \bar{\nabla}^2_{kl}s}\nonumber\\
  &-\alpha(\alpha+1) F_*^{-\alpha-2}\dot{F}_*^{pq}\frac{\partial  (\mathcal{W}_X^{-1})_{pq}}{\partial \bar{\nabla}^2_{ij}s}\dot{F}_*^{rt}\frac{\partial  (\mathcal{W}_X^{-1})_{rt}}{\partial \bar{\nabla}^2_{kl}s}
\end{align}
The uniform bound on $F$ and Lemma \ref{s5:lem-kapp} imply that there exists a constant $C>0$ such that  $ 0<C^{-1}I\leq ~(\dot{F}_*^{ij})~\leq CI$. The estimates on $s,\bar{\nabla}s$ and $F$ then imply that
 \begin{equation*}
   \lambda I\leq(\dot{G}^{ij})\leq \Lambda I
 \end{equation*}
for some constants $\lambda,\Lambda>0$. By the concavity of $F_*$ and $\alpha>0$, from \eqref{s6:dd-G} the operator $G$ is concave with respect to $\bar{\nabla}^2s$. Since we have uniform $C^2$ estimates on $s$ in space-time, we can apply the H\"{o}lder estimate of \cite{Krylov} (as in \cite{And-Wei2017,Cab-Sin2010}) to obtain the $C^{2,\alpha}$ estimate on $s$ and $C^{\alpha}$ estimate on $\partial_ts$ for some $\alpha\in (0,1)$ in space-time. By the parabolic Schauder theory \cite{lieberman1996}, we can derive estimates on all higher order derivatives of $s$. A standard continuation argument yields that the flow \eqref{flow-VMCF} exists for all time $[0,\infty)$.
\begin{prop}
Let $M_t$ be a smooth h-convex solution to the flow \eqref{flow-VMCF} with $\alpha>0$ and $\phi(t)$ given by \eqref{s1:phit}. If $f$ is inverse concave and $f_*$ approaches zero on the boundary of $\Gamma_+$, then the solution $M_t$ exists for all time $t\in[0,\infty)$.
\end{prop}

\section{Smooth convergence}\label{sec:conv}

In this section, we will use the Alexandrov reflection method to show that the solution of \eqref{flow-VMCF} converges smoothly to a geodesic sphere as $t\to\infty$.

Let $\gamma$ be a geodesic line in $\mathbb{H}^{n+1}$, and let $H_{\gamma(s)}$ be the totally geodesic hyperbolic $n$-plane in $\mathbb{H}^{n+1}$ which is perpendicular to $\gamma$ at $\gamma(s), s\in \mathbb{R}$. We use the notation $H_{s}^+$ and $H_s^-$ for the half-spaces in $\mathbb{H}^{n+1}$ determined by $H_{\gamma(s)}$:
\begin{equation*}
  H_s^+:=\bigcup_{s'\geq s}H_{\gamma(s')},\qquad  H_s^-:=\bigcup_{s'\leq s}H_{\gamma(s')}
\end{equation*}
For a bounded domain $\Omega$ in $\mathbb{H}^{n+1}$, denote
\begin{equation*}
  \Omega^+(s)=\Omega\cap H_s^+,\qquad  \Omega^-(s)=\Omega\cap H_s^-.
\end{equation*}
The reflection map across $H_{\gamma(s)}$ is denoted by $R_{\gamma,s}$. We define
\begin{equation}\label{s6:s-def}
  S_{\gamma}^+(\Omega):=\inf\{s\in \mathbb{R}~|~R_{\gamma,s}(\Omega^+(s))\subset\Omega^-(s)\}.
\end{equation}
\begin{lem}\label{s6:lem-mono-S}
For any geodesic line $\gamma$ in $\mathbb{H}^{n+1}$, $S_{\gamma}^+(\Omega_t)$ is strictly decreasing along the flow \eqref{flow-VMCF} unless $R_{\gamma,\bar s}(\Omega_t)=\Omega_t$ for some $\bar s\in\RR$
(in which case $S_{\gamma}^+(\Omega_t)=\bar s$ for all $t$).
\end{lem}

\proof
Fix $t_0\in [0,\infty)$, and set $\bar s=S_\gamma^+(\Omega_{t_0})$.   By definition we have
$R_{\gamma,\bar s}(\Omega_{t_0}^+(\bar s))\subset \Omega_{t_0}^-(\bar s)$, and since $s$ cannot be decreased below $\bar s$ we must have one of two possibilities:  Either (i) there is a point $\bar x\in R_{\gamma,\bar s}(\partial\Omega_{t_0}^+(\bar s))\cap \partial\Omega_{t_0}^-(\bar s)$ which is not in $H_{\gamma(\bar s)}$, or (ii) there is a point $\bar x$ in $\partial\Omega_{t_0}\cap H_{\gamma(\bar s)}$ such that $T_{\bar x}\partial\Omega_{t_0}$ is mapped to itself by $R_{\gamma,\bar s}$.

For case (i), since both $M_{t}^{-}(s)$ and $R_{\gamma(s)}(M_{t}^+(s))$ are strictly \emph{h-convex}, we locally express them as graphs of functions $r^-(\theta,{t})$ and $r^+(\theta,t)$ over a domain $U$ of a geodesic sphere for $t$ sufficiently close to $t_0$. Define $\omega(\theta,t):=r^-(\theta,t)-r^+(\theta,t)$. Then $\omega(\theta,t_0)$ is nonnegative for $\theta\in U$  and there exists a point $\theta_0\in U$ where the minimum $\omega(\theta_0,t_0)=0$ is achieved.  We will argue below using the strong maximum principle that $\omega$ vanishes identically, and it follows from this that $M_{t_0}^{-}(s)$ coincides with $R_{\gamma(s)}(M_{t_0}^+(s))$ and hence $M_{t_0}$ is reflection symmetric across the totally geodesic plane $H_{\gamma(s)}$ at $s=S_{\gamma}^+(\Omega_{t_0})$.

In order to apply the strong maximum principle we first recall the graphical representation of hypersurfaces in $\mathbb{H}^{n+1}$. Let $M\subset \mathbb{H}^{n+1}$ be a hypersurface which can be written as a radial graph over a sphere $S^n$, i.e., $M=\{(\theta,r(\theta)):~\theta\in S^n\}$ for a smooth function $r$ on $S^n$. Let $\{\theta^i\}, i=1,\cdots,n$ be a local coordinate system on $S^n$. The induced metric on $M$ from $\mathbb{H}^{n+1}$ takes the form
\begin{equation*}
  g_{ij}~=~D_irD_jr+\sinh^2r\sigma_{ij},
\end{equation*}
where $\sigma_{ij}$ denotes the standard metric on $S^n$. Denote
 \begin{equation}\label{s6:v}
  v=\sqrt{1+{|Dr|^2}/{\sinh^2r}},\quad \mathrm{and }\quad  \tilde{\sigma}^{jk}~=~\sigma^{jk}-\frac{D^jrD^kr}{v^2\sinh^2r}
\end{equation}
where $\sigma^{jk}$ is the inverse matrix of $\sigma_{jk}$ and $D^jr=\sigma^{jk}D_kr$. Then the Weingarten matrix $(h_i^j)$ can be expressed as (cf. \cite{Gerh2006})
\begin{equation}\label{s6:h}
  h_i^j=~\frac{\coth r}{v}\delta_i^j+\frac{\coth r}{v^3\sinh^2 r}D^jrD_ir-\frac {\tilde{\sigma}^{jk}}{v\sinh^2r}D_{k}D_ir.
\end{equation}
Up to a tangential diffeomorphism, the flow equation \eqref{flow-VMCF} is equivalent to the following scalar parabolic PDE
\begin{equation}\label{s6:graph-flow}
  \frac{\partial r}{\partial t}=~(\phi(t)-F^{\alpha}(h_i^j))\sqrt{1+{|Dr|^2}/{\sinh^2r}}.
\end{equation}
Denote the right hand side of \eqref{s6:graph-flow} by $ \Phi(D^2r, Dr, r,\theta,t)$.

We now come back to the $M_{t}^{-}(s)$ and $R_{\gamma(s)}(M_{t}^+(s))$ and $\omega(\theta,t)=r^-(\theta,t)-r^+(\theta,t)$. Since the flow \eqref{flow-VMCF} is invariant under reflection, by \eqref{s6:graph-flow} the function $\omega(\theta,t)$ satisfies the following equation
\begin{align}\label{s6:evl-omega}
  \frac{\partial\omega}{\partial t}=&~  \Phi(D^2r^-, Dr^-, r^-,\theta,t)-  \Phi(D^2r^+, Dr^+, r^+,\theta,t)\nonumber\\
  =&~\frac{\partial\Phi}{\partial D^2_{ij}r}(\chi,\xi,\eta)D^2_{ij}\omega+\frac{\partial\Phi}{\partial D_{i}r}(\chi,\xi,\eta)D_{i}\omega+\frac{\partial\Phi}{\partial r}(\chi,\xi,\eta)\omega,
\end{align}
where $(\chi,\xi,\eta)=(\rho D^2r^-+(1-\rho)D^2r^+, \rho Dr^-+(1-\rho)Dr^+, \rho r^-+(1-\rho)r^+)$ for some $\rho\in[0,1]$. The uniform estimates in \S \ref{sec:LTE} implies that the equation \eqref{s6:evl-omega} is uniformly parabolic, i.e., there exist constants $\lambda,\Lambda>0$ such that
\begin{equation*}
  \lambda I\leq (\frac{\partial\Phi}{\partial D^2_{ij}r})\leq \Lambda I.
\end{equation*}
The coefficients $\frac{\partial\Phi}{\partial D_{i}r}$ and $\frac{\partial\Phi}{\partial r}$ are uniformly bounded and smooth. Since $\omega(\theta,t_0)$ is nonnegative and is positive somewhere in $U$, the strong maximum principle applied to the equation \eqref{s6:evl-omega} yields that $\omega>0$ everywhere in $U$ for $t>t_0$, unless it is identically zero. This in turn implies that $S_{\gamma}^+(\Omega_t)<S_{\gamma}^+(\Omega_{t_0})$ for $t>t_0$.

The discussion for case (ii) is similar. We again write $M_{t}^{-}(s)$ and $R_{\gamma(s)}(M_{t}^+(s))$ locally as graphs of functions $r^-(\theta,{t})$ and $r^+(\theta,t)$ over a domain $U$ of a geodesic sphere for $t$ sufficiently close to $t_0$, and then apply the boundary strong maximum principle (the Hopf Lemma).
\endproof

Now we prove that the flow \eqref{flow-VMCF} converges smoothly to a geodesic sphere. For any fixed $\tau$, we define the flow $\Omega_{\tau}(t)$ by $\Omega_{\tau}(t):=\Omega_{t+\tau}$. The uniform estimates in \S \ref{sec:LTE} imply that there exists a sequence of $\tau_k\to\infty$ such that the families $\Omega_{\tau_k}(t)$ converge smoothly to a limiting flow $\Omega_{\infty}(t), t\in[0,\infty)$ which is again a solution of \eqref{flow-VMCF}. By the monotonicity of $S_{\gamma}^+(\Omega_t)$ proved in Lemma \ref{s6:lem-mono-S}, we have
\begin{equation}\label{s6:limf-1}
  S^+_{\gamma}(\Omega_{\infty}(t))=\lim_{t'\to\infty}S_{\gamma}^+(\Omega_{t'})
\end{equation}
which exists by the monotonicity. The right hand side of \eqref{s6:limf-1} is independent of $t$ and is finite. We conclude that the limiting flow $\Omega_{\infty}(t)$ is symmetric under reflection across a totally geodesic hyperplane $H_{\gamma}$ which is perpendicular to the geodesic line $\gamma$.  Since $\gamma$ is arbitrary, we conclude that $\Omega_{\infty}(t)$ is a geodesic sphere. This implies the subsequential smooth convergence of $\Omega_t$ to a geodesic sphere of radius $r_{\infty}=f_k^{-1}(W_k(\Omega_0))$.

The linearisation of the flow \eqref{flow-VMCF} about the geodesic sphere of radius $r_{\infty}$ can be used to deduce the stronger convergence results. The hypersurface near the geodesic sphere can be written as a graph of a smooth function $r$ over the geodesic sphere. Setting $r=r_{\infty}(1+\epsilon\eta)$. The linearised equation of the flow \eqref{s6:graph-flow} about the geodesic sphere of radius $r_{\infty}$ is given by
\begin{equation}\label{s6:graph-Linear}
  \frac{\partial\eta}{\partial t}=~\frac{\alpha\coth^{\alpha-1}r_{\infty}}{n\sinh^2r_{\infty}}\left(\Delta\eta+n\eta-\frac n{|S^n|}\int_{S^n}\eta d\sigma\right)
\end{equation}
This equation is the same as that for mixed volume preserving mean curvature flow in \cite[Eqn.(21)]{McC2004} (see also \cite[\S 12]{And2001-Aniso}). Thus the same argument as in \cite{And2001-Aniso} gives that the flow \eqref{flow-VMCF} converges exponentially in smooth topology to the geodesic sphere with radius $r_{\infty}$. This completes the proof of Theorem \ref{thm1-1}.

\appendix
\section{Smooth convergence: $f=E_k^{1/k}$}
In this appendix, we provide an alternative approach to the proof of smooth convergence to a geodesic sphere for the volume preserving flow \eqref{flow-VMCF} with $f=E_k^{1/k}$ and $\alpha>0$. The key ingredient is the monotonicity of the quermassintegral $W_k(\Omega_t)$.

\begin{lem}\label{lem-monot}
For any integer $k\in \{1,\cdots,n\}$, let $M_t$ be a smooth convex solution of the volume preserving flow \eqref{flow-VMCF} with $f=E_k^{1/k}$ and $\alpha>0$. Denote $\Omega_t$ the domain enclosed by $M_t$. Then $W_k(\Omega_t)$ is monotone decreasing in time $t$.
\end{lem}
\proof
The function $\phi(t)$ in \eqref{flow-VMCF} for the volume preserving flow is defined as in \eqref{s2:phit}. From the evolution equation \eqref{s2:evl-Wk-1} for the quermassintegral of $\Omega_t$:
\begin{equation}\label{s4:evl-Wk-1}
  \frac d{dt}W_k(\Omega_t)=~\frac {n+1-k}{n+1}\int_{M_t}E_k\left(\phi(t)-F^{\alpha}\right)d\mu_t
\end{equation}
and $F=E_{k}^{1/k}$, we have
\begin{equation*}
  \frac d{dt}W_k(\Omega_t)=~\frac {n+1-k}{n+1}\left(\frac 1{|M_t|}\int_{M_t}E_k\int_{M_t}E_k^{\alpha/k}-\int_{M_t}E_k^{1+\alpha/k}d\mu_t\right).
\end{equation*}
Then the monotonicity of $W_k(\Omega_t)$ follows immediately from the Jensen's inequality
\begin{equation}\label{s3:Jensen}
  \int_{M_t}E_k^{(\alpha+k)/k}d\mu_t~\geq ~\frac 1{|M_t|}\int_{M_t}E_{k}d\mu_t\int_{M_t}E_k^{\alpha/k}d\mu_t.
\end{equation}
\endproof
If the initial hypersurface $M_0$ is \emph{h-convex}, then the long time existence of the flow \eqref{flow-VMCF} has been proved in \S \ref{sec:LTE}. To show the smooth convergence to a geodesic sphere, we need the following lemma. Denote
\begin{equation*}
  \bar{E}_k=~\frac 1{|M_t|}\int_{M_t}E_kd\mu_t.
\end{equation*}
\begin{lem}\label{s6:lem2}
For any integer $k\in \{1,\cdots,n\}$, let $M_t$ be a smooth $h$-convex solution of the volume preserving flow \eqref{flow-VMCF} with $f=E_k^{1/k}$ and $\alpha>0$. Then there exists a sequence of times $t_i\to\infty$ such that
\begin{align}\label{s6:1}
 \int_{M_{t_i}}\left(E_k-\bar{E}_k\right)^2d\mu_{t_i}~\to&~0, \quad\mathrm{ as}~i\to\infty
\end{align}
\end{lem}
\proof
Since $M_t$ is \emph{h-convex}, the Alexandrov-Fenchel inequality  \eqref{s2:AF1} implies that ,
\begin{equation*}
   W_k(\Omega_t)~\geq~f_k\circ f_0^{-1}(W_0(\Omega_t))=~f_k\circ f_0^{-1}(W_0(\Omega_0))~>~0,
\end{equation*}
where we used the condition $|\Omega_t|=|\Omega_0|$. By Lemma \ref{lem-monot},  $W_k(\Omega_t)$ is monotone decreasing. Then we can find a sequence of times $t_i\to\infty$ such that
\begin{equation*}
  \frac d{dt}\bigg|_{t_i}W_k(\Omega_t)~\to~0,\quad\mathrm{ as}~i\to\infty,
\end{equation*}
Then from the proof of Lemma \ref{lem-monot}, equality is attained in Jensen's inequality \eqref{s3:Jensen} as $t_i\to\infty$, or equivalently
\begin{align}\label{5-1}
  0\leq &~\int_{M_{t_i}}E_k^{\alpha/k}\left(E_k-\bar{E}_k\right)d\mu_{t_i}\nonumber\\
  =&\int_{M_{t_i}}\left(E_k^{\alpha/k}-\bar{E}_k^{\alpha/k}\right)\left(E_k-\bar{E}_k\right)d\mu_{t_i}~\to~0, \quad\mathrm{ as}~i\to\infty.
\end{align}
Since
\begin{equation}\label{5-2}
  \left(E_k^{\alpha/k}-\bar{E}_k^{\alpha/k}\right)\left(E_k-\bar{E}_k\right)~\geq~C\left(E_k-\bar{E}_k\right)^2,
\end{equation}
as in the proof of Lemma 6.1 of \cite{And-Wei2017}, then \eqref{s6:1} follows from \eqref{5-1} and \eqref{5-2} immediately.
\endproof

The uniform estimates on all the higher derivatives of the Weingarten map $(h_i^j)$ implies that for any sequence of times $t_i\to\infty$, there exists a subsequence (still denoted by $t_i$) such that $M_{t_i}$ converges to a limit hypersurface $M_{\infty}$ (up to an ambient isometry). Thus by Lemma \ref{s6:lem2}, we can find a sequence of times $t_i\to\infty$ such that $M_{t_i}$ converges smoothly to a geodesic sphere up to an isometry. For any other sequence of times $t_j\to\infty$, we can also find a subsequence of time (labeled by the same $t_j$) such that $M_{t_j}$ converges to a limit $\hat{M}_{\infty}=\partial\hat{\Omega}_{\infty}$. The monotonicity of  $W_k(\Omega_t)$ then yields $W_k(\hat{\Omega}_{\infty})=W_k(B^{n+1})$, which  implies that $\hat{\Omega}_{\infty}$ is a geodesic ball by the equality case of Theorem \ref{s2:thm1}. Since the above argument works for any sequence, we conclude that the whole family of $M_t$ converges to a geodesic sphere as $t\to\infty$ up to an isometry. The exponential convergence is the same as in \S \ref{sec:conv}.

\bibliographystyle{amsplain}

\end{document}